# OPTIMAL DESIGNS WHICH ARE EFFICIENT FOR LACK OF FIT TESTS

By Wolfgang Bischoff and Frank Miller

*Catholic University of Eichstätt–Ingolstadt and AstraZeneca*

Linear regression models are among the models most used in practice, although the practitioners are often not sure whether their assumed linear regression model is at least approximately true. In such situations, only designs for which the linear model can be checked are accepted in practice. For important linear regression models such as polynomial regression, optimal designs do not have this property. To get practically attractive designs, we suggest the following strategy. One part of the design points is used to allow one to carry out a lack of fit test with good power for practically interesting alternatives. The rest of the design points are determined in such a way that the whole design is optimal for inference on the unknown parameter in case the lack of fit test does not reject the linear regression model.

To solve this problem, we introduce efficient lack of fit designs. Then we explicitly determine the $\mathbf{e}_k$-optimal design in the class of efficient lack of fit designs for polynomial regression of degree $k-1$.

**1. Introduction.** Linear regression models are among the models most used in practice. Such a parametric assumption for the regression function is very attractive among practitioners, although they are often not sure whether their assumed linear regression model is at least approximately true. Therefore, if a design can be chosen (according to which the data are sampled), the practitioners spread out the design points over the whole experimental region. For important linear regression models such as polynomial regression, such designs and classical optimal designs are quite different. Even more serious when using such an optimal design, deviations from the assumed polynomial regression model are not detectable. In this paper we address these concerns.

To explain the above mentioned problem in more detail, let us consider the linear regression model

$$\mathbf{Y} = X\boldsymbol{\theta} + \boldsymbol{\varepsilon}, \tag{1.1}$$

where $\mathbf{Y} = (Y_1, \ldots, Y_n)^\top$ is the vector of observations, $X$ is the design (model) matrix, $\boldsymbol{\theta} \in \mathbb{R}^k$ is an unknown parameter vector and $\boldsymbol{\varepsilon} = (\varepsilon_1, \ldots, \varepsilon_n)^\top$ is the vector of errors. In this paper we assume that $\varepsilon_1, \ldots, \varepsilon_n$ are independent, identically distributed real random variables with expectation 0 and unknown variance $\sigma^2 \in (0, \infty)$. Furthermore, we assume $X = X_n = (f(x_1), \ldots, f(x_n))^\top$, where $(x_1, \ldots, x_n) \in \mathcal{E}^n$ is a design for $n$ observations, $\mathcal{E} = [a, b] \subseteq \mathbb{R}$ is some compact interval and $f = (f_1, \ldots, f_k)^\top : \mathcal{E} \to \mathbb{R}^k$ is a vector of known continuous regression functions that have bounded variation.

It is common to solve classical experimental design problems for linear regression models of the form (1.1) in an approximate way. To do this, one identifies an arbitrary design with $n$ design points $(x_1, \ldots, x_n) \subseteq \mathcal{E}^n$ with the probability measure $\xi_n := \frac{1}{n} \sum_{i=1}^n \delta_{x_i}$, where $\delta_t$ is the Dirac measure in $t \in \mathcal{E}$. Then it is more feasible to solve a classical design problem in the set of all probability measures on $\mathcal{E}$ instead in $\Xi_n := \{\frac{1}{n} \sum_{i=1}^n \delta_{x_i} | (x_1, \ldots, x_n) \in \mathcal{E}^n\}$, or in $\bigcup_{n=1}^\infty \Xi_n$. If we are interested in inference on the parameter vector $K^\top \boldsymbol{\theta}$, where $K \in \mathbb{R}^{k \times s}$, $1 \leq s \leq k$, $\text{rank}(K) = s$, then, given the design $\xi$, the variance/covariance matrix of the best linear unbiased (least squares) estimator for $K^\top \boldsymbol{\theta}$ is given by

$$\frac{\sigma^2}{n} K^\top M(\xi)^{-1} K, \qquad M(\xi) := \int_\mathcal{E} f f^\top \, d\xi$$

[provided that $M(\xi)$ is invertible, otherwise we have to deal with generalized inverses]; see [17], page 65. Hence, the "quality" of statistical inference on the unknown parameter $\boldsymbol{\theta}$ depends on the choice of the design $\xi$. The most interesting optimality criteria, such as the $\phi_p$-criteria, are functions of the information matrix $C(\xi) := (K^\top M(\xi)^{-1} K)^{-1}$. By a famous theorem of Carathéodory (see, e.g., [22]), the optimal designs for a criterion that is a function of the information matrix can be determined in the set of probability measures with $k(k+1)/2$ mass points instead of in the set of all probability measures. Indeed, the classical optimal designs for polynomial regression were determined in this set of probability measures just mentioned. Therefore, in case the number $k$ of different regression functions is small (which is most interesting in practice), the optimal designs have only a small number of different design (mass) points. For instance, for polynomial regression of degree $k-1$, the known classical $D$-, $G$-, $A$- and $\mathbf{e}_k$-optimal [with $\mathbf{e}_k = (0, \ldots, 0, 1)^\top \in \mathbb{R}^k$] designs have only $k$ different design points. There is a huge literature on classical optimal designs in case model (1.1) is true; see, for example, [9, 12, 14, 17, 22].

To check models, Wiens [24] introduced optimal lack of fit- (LOF-) designs based on the power of a given lack of fit test. Wiens [24] considered the usual lack of fit test for regression models of the form (1.1). This approach was recently generalized by Biedermann and Dette [1].

It is of practical interest to develop (in some sense) optimal designs that offer the possibility to check whether the assumed model is true (or is at least not so far away from the true one) and with which one makes good inference on $\boldsymbol{\theta}$ in case the assumed regression model is true. To establish designs which fulfill the above two postulates, we combine both ideas. For that, let $g:\mathcal{E} \to \mathbb{R}$ be the true but unknown regression function. Then the hypothesis

$$(1.2) \qquad H_0: \exists\, \boldsymbol{\theta} = (\theta_1, \ldots, \theta_k)^\top \in \mathbb{R}^k \qquad \text{with } g = f^\top \boldsymbol{\theta}$$

is of practical main interest. A test for this hypothesis is called "lack of fit"-test (LOF-test). Loosely speaking, our approach is then as follows. Let $r \in [0,1]$ and let a certain lack of fit test be given. Then an optimal LOF-design has the "best" power for certain alternatives (of regression functions) that are separated from $H_0$; see Section 2.2. At first, we determine the class of $r$-efficient LOF-designs, that is, all designs for which the efficiency with respect to the optimal LOF-design is $r$. Hence, for $r = 1$, we obtain the optimal LOF-designs. Note that the known classical optimal designs for polynomial regression models, however, often belong to the class of 0-efficient LOF-designs. Next, in the class of LOF-designs which are at least $r$-efficient, the optimal design with respect to statistical inference on $\boldsymbol{\theta}$ is determined under the assumption that the model (1.1) is true. In practice, the value $r$ of efficiency has to be chosen by the practitioner according to the statistical problem under consideration. Then the data are sampled according to the $r$-efficient-optimal LOF-design. In case model (1.1) is rejected by the LOF-test, the data are useless for inference on $\boldsymbol{\theta}$. In case model (1.1) is not rejected, the data are used for inference on $\boldsymbol{\theta}$.

In Section 2 we compare ours with other approaches that handle the problem under consideration. Furthermore, the class of $r$-efficient LOF-designs is established there. Then we determine $\mathbf{e}_k$-optimal designs in the class of at least $r$-efficient LOF-designs for polynomial regression of degree $k-1$ in Section 3. For corresponding results for a general linear regression model (1.1), see [15]. Moreover, we show the relation to Bayesian optimal designs. Proofs are postponed to the Appendix.

## 2. LOF-designs.

2.1. *LOF-tests and an asymptotic aspect of LOF-designs.* Let us again consider the assumed linear model (1.1), $\mathbf{Y} = X\boldsymbol{\theta} + \boldsymbol{\varepsilon} = (f(x_{n1}), \ldots,$ $f(x_{nn}))^\top \boldsymbol{\theta} + \boldsymbol{\varepsilon}$, the true model, $\mathbf{Y} = (g(x_{n1}), \ldots, g(x_{nn}))^\top + \boldsymbol{\varepsilon}$, and the hypothesis (1.2), $H_0 : \exists \boldsymbol{\theta} = (\theta_1, \ldots, \theta_k)^\top \in \mathbb{R}^k$ with $g = f^\top \boldsymbol{\theta}$. Many LOF-tests are possible for $H_0$. Wiens [24] considered the usual LOF-test that has a noncentral $F$-distribution if the errors are normally distributed. This holds asymptotically as well under mild conditions; see [24]. If the regression function and the density of the design $\xi$ fulfill some smoothness assumptions, then tests based on nonparametric estimation of the regression function are possible. Biedermann and Dette [1] proposed three LOF-tests of this kind. Every test mentioned above has an asymptotic (for $\xi_n$ converging in some sense to $\xi$, as $n \to \infty$) power, which is a function of

$$B(g,\xi) := \frac{1}{\sigma^2} \int_{\mathcal{E}} ((\mathrm{pr}^{L^2(\xi)}_{[f_1, \ldots, f_k]^\perp} g)(x))^2 \xi(dx),$$

where $\mathrm{pr}^{L^2(\xi)}_{[f_1, \ldots, f_k]^\perp}$ is the orthogonal projector onto $[f_1, \ldots, f_k]^\perp$ in $L^2(\xi)$; see [1, 24]. The greater $B(g,\xi)$ is, the greater is the asymptotic power. It is easy to argue that we have to put $B(g,\xi) := \infty$ for $g \notin L^2(\xi)$. In practice, however, it is sufficient to consider regression functions that have bounded variation. Additionally, we assume that the regression function is regular, which here means $g$ is continuous from the left in $b$ and continuous from the right in $[a, b)$. We denote this class of functions by $BV(\mathcal{E}) = BV[a, b]$. In the sequel we assume that the true regression function $g$ is an element of $BV(\mathcal{E})$.

Opposed to the approximation approach for the classical design problem, an asymptotic statistical argument is decisive for the set of probability measures in which the LOF-design problem should be solved. Bischoff [2, 3] and Bischoff and Miller [4] consider partial sum processes of regression models. These papers imply that a regression function $g \in BV(\mathcal{E})$ which does not belong to the assumed regression model can be detected asymptotically if the sequence of exact designs converges uniformly to an asymptotic design (probability measure) that has an absolutely continuous part with positive density which belongs to $BV(\mathcal{E})$. Uniform convergence of probability distributions means uniform convergence of the corresponding distribution functions. Therefore, as a class of interesting designs, we consider the set $\Xi$ of probability measures on $\mathcal{E}$ that can be decomposed into a finitely discrete part and an absolutely continuous part with respect to the Lebesgue measure whose density belongs to $BV(\mathcal{E})$. Moreover, the papers mentioned above additionally show that $B(g,\xi)$ is a suitable measure for how well an alternative regression function can be detected. Indeed, $\sqrt{B(g,\xi)}$ is the norm of the residual sum limit process with respect to the reproducing kernel Hilbert space of the corresponding limit process. Before we continue, we have to approximate $\xi \in \Xi$ by a design in $\Xi_n$ if $n$ is the fixed number of observations. To this end, let $F_0(t) := \xi((-\infty, t]), t \in \mathcal{E}$, be the distribution function of $\xi$, and let $Q_0$ be the right continuous inverse of $F_0$. Then the design $\xi_n = \frac{1}{n}\sum_{i=1}^n \delta_{x_{ni}} \in \Xi_n$ with

$$(2.1) \qquad x_{ni+1} := Q_0(i/(n-1)), \qquad i = 0,\ldots, n-1,$$

has the property that $\xi_n$ converges uniformly to $\xi$. These designs correspond to designs defined by Sacks and Ylvisaker [20].

There are several other approaches that handle the problem under consideration. For polynomial regression of fixed degree, for instance, Dette [8] and Pukelsheim and Rosenberger [18] consider designs for a polynomial of higher degree (as an alternative) and a mixture of two optimality criteria. Box and Draper [7] (see also [11]) look for designs minimizing a bias. Montepiedra and Yeh [16] uses a sequential approach. Biswas and Chaudhuri [6] select the correct model from a known finite family of nested linear models and estimate the parameters associated with that model. But the optimal designs of these approaches do not have the property that each fixed alternative of the class mentioned above can be discovered as $n \to \infty$.

2.2. *Optimal LOF-designs.* For the definition of LOF-efficiency, it is technically necessary to consider a subset of alternatives which is separated from the hypothesis $H_0$. Following Biedermann and Dette [1], let $v \in BV(\mathcal{E})$ be a weight function which gives more weight to those points of $\mathcal{E}$ for which a deviation is more serious, and let $\lambda$ be Lebesgue measure. Then we consider $\{h \in BV(\mathcal{E})| \int_\mathcal{E} h^2 v\, d\lambda \geq c,\ \int_\mathcal{E} fhv\, d\lambda = \mathbf{0}_k\}$, $c > 0$ fixed, as a set of alternatives which is separated from the hypothesis $H_0$. Let $\tilde{\lambda}$ be the uniform distribution on $\mathcal{E}$, that is, $\tilde{\lambda} = (\lambda(\mathcal{E}))^{-1} \cdot \lambda|_\mathcal{E}$. By choosing $v : \mathcal{E} \to [0, \infty)$ with

$$(2.2) \qquad \int_\mathcal{E} v\, d\tilde{\lambda} = 1$$

and $c > 0$ suitably, the above set of alternatives can be written as

$$\mathcal{F}_{v,c} := \left\{ h \in BV(\mathcal{E}) \Big| \int_\mathcal{E} h^2 v\, d\tilde{\lambda} \geq c, \int_\mathcal{E} fhv\, d\tilde{\lambda} = \mathbf{0}_k \right\}.$$

In the sequel let $\mathcal{F}_{v,c}$ be fixed and let $v$ satisfy (2.2). Following Wiens [24], we continue with a maximin approach.

DEFINITION 2.1. (a) [24]. A design $\xi_0 \in \Xi$ is called LOF-optimal if

$$\max_{\xi \in \Xi} \min_{h \in \mathcal{F}_{v,c}} B(h, \xi) = \min_{h \in \mathcal{F}_{v,c}} B(h, \xi_0).$$

(b) The relative LOF-efficiency of a design $\xi_1 \in \Xi$ is

$$\mathrm{Eff}_{\mathrm{LOF}}(\xi_1) = \min_{h \in \mathcal{F}_{v,c}} B(h, \xi_1) \Big/ \min_{h \in \mathcal{F}_{v,c}} B(h, \xi_0) \in [0, 1],$$

where $\xi_0$ is an optimal LOF-design.

2.3. *Efficient LOF-designs.* Wiens [24] computed optimal LOF-designs for $v \equiv const$. Biedermann and Dette [1] generalized this result to arbitrary $v$. These papers imply that the optimal LOF-design is the probability measure $v \cdot \tilde{\lambda}$, where for a measure $\eta$ defined on the Borel field $\mathcal{B}$ and a Borel-measurable function $w : \mathbb{R} \to [0, \infty)$ the measure $w \cdot \eta$ is defined by $(w \cdot \eta)(A) := \int_A w \, d\eta$, $A \in \mathcal{B}$. Next we give a generalization of this result. We use for two measures $\mu_1, \mu_2$ on $\mathcal{B}$ the notation $\mu_1 \leq \mu_2 \Leftrightarrow \forall B \in \mathcal{B} : \mu_1(B) \leq \mu_2(B)$.

THEOREM 2.2. *The set of designs with relative LOF-efficiency greater than or equal to $r \in [0,1]$ is given by $\Upsilon_v[r] := \{\xi \in \Xi | rv \cdot \tilde{\lambda} \leq \xi\}$.*

Given $\xi \notin \Upsilon_v[r]$, the main part of the proof is to construct a regression function $h_0 \in \mathcal{F}_{v,c}$ with $(\mathrm{Eff}_{\mathrm{LOF}}(\xi) \leq) B(h_0, \xi) / \min_{h \in \mathcal{F}_{v,c}} B(h, v \cdot \tilde{\lambda}) < r$. The proof is mainly along the lines of Wiens [24] and Biedermann and Dette [1] and therefore is omitted. An immediate consequence of Theorem 2.2 is the following interesting corollary.

COROLLARY 2.3. *Let $\xi \in \Xi$ and $r := \sup\{t | tv \cdot \tilde{\lambda} \leq \xi\} \in [0,1]$. Then the LOF-efficiency of $\xi$ is equal to $r$.*

**3. $\mathbf{e}_k$-optimal designs in $\Upsilon_v[r]$ for polynomial regression.** In this section we calculate $\mathbf{e}_k$-optimal designs in $\Upsilon_v[r]$, where $v$ is an arbitrary weight function, in case the parametric model (1.1) is the polynomial regression model of degree $k-1$, that is, $f(x) = (1, x, \ldots, x^{k-1})^\top$, $k \geq 2$, where $\mathbf{e}_k = (0, \ldots, 0, 1)^\top \in \mathbb{R}^k$. As discussed in [15], we can consider the experimental region $\mathcal{E} = [-1, 1]$ for $\mathbf{e}_k$-optimality, without loss of generality. Furthermore, the $\mathbf{e}_k$-optimal design in $\Upsilon_v[r]$ for polynomial regression is unique; see also [15]. The proof of the following theorem can be found in the Appendix.

THEOREM 3.1. *A design $\xi \in \Upsilon_v[r]$ is optimal for $\mathbf{e}_k^\top \boldsymbol{\theta}$ in $\Upsilon_v[r]$ if and only if*

$$
\begin{aligned}
&\forall y \in [-1, 1], \ \forall z \in \{z \in [-1, 1] | \xi(\{z\}) > 0\} : \\
&(\mathbf{e}_k^\top M(\xi)^{-1} f(y))^2 \leq (\mathbf{e}_k^\top M(\xi)^{-1} f(z))^2.
\end{aligned}
\tag{3.1}
$$

*Moreover, the optimal design for $\mathbf{e}_k^\top \boldsymbol{\theta}$ in $\Upsilon_v[r]$ has the form $rv \cdot \tilde{\lambda} + (1-r) \times \sum_{i=1}^\ell p_i \delta_{t_i}$, where $t_1, \ldots, t_\ell \in [-1, 1]$ are $\ell(\leq k)$ different points with $|\mathbf{e}_k^\top \times M(\xi)^{-1} f(t_i)| = \max_{y \in [-1,1]} \mathbf{e}_k^\top M(\xi)^{-1} f(y)$ and $p_1, \ldots, p_\ell \in (0, 1]$ are suitable values with $\sum_{i=1}^\ell p_i = 1$.*

Note that our design problem is related to a Bayesian design problem. For that, let $rv \cdot \tilde{\lambda} + (1-r)\zeta^*$ be the $\mathbf{e}_k^\top \boldsymbol{\theta}$-optimal design in $\Upsilon_v[r]$ and let $M_0 := M(v \cdot \tilde{\lambda})$. Then $\zeta^*$ is the Bayesian optimal design for $\mathbf{e}_k^\top \boldsymbol{\theta}$ of the Bayesian design problem with a priori information $M_0$. This means that $\zeta^*$ maximizes the information matrix' $rM_0 + (1-r)M(\zeta)$ with respect to $\zeta \in \Xi$; see [10], Section 5, [8], [17], Chapter 11, and [15].

Kiefer and Wolfowitz [13] showed that $\frac{1}{2(k-1)}(\delta_{-1} + \delta_1) + \frac{1}{k-1} \times \sum_{i=1}^{k-2} \delta_{\cos(\pi(k-1-i)/(k-1))}$ is an $\mathbf{e}_k$-optimal design in $\Xi$; see also [23]. The support points of this optimal design are the extremal points in $[-1, 1]$ of the Chebyshev polynomial $T_{k-1}$ of degree $k-1$. See [19] and [21] for properties of Chebyshev polynomials $T_n(x), n \in \mathbb{N}_0, x \in [-1, 1]$. The proof of the following main result can be found in the Appendix.

THEOREM 3.2. *Let $\nu_i = 1/2$ if $i \in \{0, k-1\}$, $\nu_i = 1$ if $i \in \{1, \ldots, k-2\}$, let $p_i = \frac{\nu_i}{k-1} - rq_i$, $i = 0, \ldots, k-1$, where*

$$q_i = \frac{\nu_i}{2(k-1)} \sum_{j=0}^{2k-3} \cos\left(\frac{j(k-1-i)\pi}{k-1}\right) \int_{-1}^{1} T_j(x) v(x) \, dx,$$

*and let*

(3.2) $\quad \alpha_0^{(k)} := \min\{\nu_i/((k-1)q_i) | i = 0, \ldots, k-1 \text{ with } q_i > 0\}.$

*Then for $r \in [0, \alpha_0^{(k)}]$, the design $\xi = rv \cdot \tilde{\lambda} + \sum_{i=0}^{k-1} p_i \delta_{\cos(\pi(k-1-i)/(k-1))}$ is $\mathbf{e}_k$-optimal in $\Upsilon_v[r]$.*

Next let us consider this result for two special weight functions.

EXAMPLE 3.3. For the weight function $v(x) = 2/(\pi\sqrt{1-x^2})$, Theorem 3.2 can be simplified substantially. Since

$$\int_{-1}^{1} T_j(x)v(x) \, dx = \frac{2}{\pi} \int_{-1}^{1} T_j(x) \frac{dx}{\sqrt{1-x^2}} = \begin{cases} 2, & \text{if } j = 0, \\ 0, & \text{otherwise} \end{cases}$$

(see [19], page 35), we get $p_i = (1-r)\frac{\nu_i}{k-1}$ and $\alpha_0^{(k)} = 1$, where $\nu_i = 1/2, i \in \{0, k-1\}$, $\nu_i = 1, i \in \{1, \ldots, k-2\}$. Hence, for arbitrary $r \in [0, 1]$, the $\mathbf{e}_k$-optimal design in $\Upsilon_v[r]$ is

$$\xi = rv \cdot \tilde{\lambda} + \frac{1-r}{k-1}\left(\frac{1}{2}(\delta_{-1} + \delta_1) + \sum_{i=1}^{k-2} \delta_{\cos(\pi(k-1-i)/(k-1))}\right).$$

EXAMPLE 3.4. Finally, we specialize Theorem 3.2 for $v \equiv 1$ by using standard results for Chebyshev polynomials. Let $\nu_i = 1/2$ for $i \in \{0, k-1\}$, let $\nu_i = 1$ for $i \in \{1, \ldots, k-2\}$, let

$$p_i = p_{k-1-i} = \frac{\nu_i}{k-1} - rq_i$$

$$= \frac{\nu_i}{k-1} - r\frac{1}{2}\frac{\nu_i}{k-1} \sum_{j=0}^{k-2} \cos\left(\frac{2ji\pi}{k-1}\right) \cdot \left(\frac{1}{2j+1} - \frac{1}{2j-1}\right)$$

for $i = 0, \ldots, k-1$, and let $\alpha_0^{(k)} = \min\{\nu_i/((k-1)q_i)|i = 0, \ldots, k-1$ with $q_i > 0\}$. Then for $r \in [0, \alpha_0^{(k)}]$ arbitrarily fixed, the design $\xi = r\tilde{\lambda} + \sum_{i=0}^{k-1} p_i \times \delta_{\cos(\pi(k-1-i)/(k-1))}$ is $\mathbf{e}_k$-optimal in $\Upsilon_v[r]$. In Table 1 we state the values $\alpha_0^{(k)} = \min\{\nu_i/((k-1)q_i)|i = 0, \ldots, k-1$ with $q_i > 0\}$ for $k = 2, 3, \ldots, 8$. Next, we state the optimal designs in $\Upsilon[r] := \Upsilon_v[r]$ for polynomial regression of certain degrees.

(a) Straight-line regression, that is, $k = 2$. Let $r \in [0,1]$ be arbitrarily fixed. Then the design $\xi = r\tilde{\lambda} + (1-r)\frac{1}{2}(\delta_{-1} + \delta_1)$ is $\mathbf{e}_2$-optimal in $\Upsilon[r]$.

(b) Quadratic regression, that is, $k = 3$. Let $r \in [0, 3/4]$ be arbitrarily fixed. Then the design $r\tilde{\lambda} + (\frac{1}{4} - \frac{r}{6})(\delta_{-1} + \delta_1) + (\frac{1}{2} - \frac{2r}{3})\delta_0$ is $\mathbf{e}_3$-optimal in $\Upsilon[r]$.

(c) Cubic regression, that is, $k = 4$. Let $r \in [0, 5/6]$ be arbitrarily fixed. Then the design $r\tilde{\lambda} + (\frac{1}{6} - \frac{r}{10})(\delta_{-1} + \delta_1) + (\frac{1}{3} - \frac{2r}{5})(\delta_{-1/2} + \delta_{1/2})$ is $\mathbf{e}_4$-optimal in $\Upsilon[r]$.

(d) Polynomial regression of degree 4, that is, $k = 5$. Let $r \in [0, 105/136]$ be arbitrarily fixed. Then the design $\xi = r\tilde{\lambda} + (\frac{1}{8} - \frac{r}{14})(\delta_{-1} + \delta_1) + (\frac{1}{4} - \frac{4r}{15})(\delta_{-1/\sqrt{2}} + \delta_{1/\sqrt{2}}) + (\frac{1}{4} - \frac{34r}{105})\delta_0$ is $\mathbf{e}_5$-optimal in $\Upsilon[r]$.

It is worth mentioning that in [5] the case that $r$ is near 1 is considered. For this case the $\mathbf{e}_k$-optimal designs in $\Upsilon[r]$ can also be explicitly calculated for $r \in [\alpha_1^{(k)}, 1]$, where $\alpha_1^{(k)}$ is some bound that can be calculated. By this result, for each $r \in [0,1]$ the $\mathbf{e}_k$-optimal designs in $\Upsilon[r]$ can be calculated for $k = 3, 4$.

## APPENDIX: PROOFS OF THEOREMS 3.1 AND 3.2

We first state an equivalence theorem which is useful for explicitly calculating optimal designs in $\Upsilon_v[r]$. Since $\{M(\xi)|\xi \in \Upsilon_v[r]\}$ is convex and compact, the proof of the following equivalence theorem is related to equivalence theorems stated in the literature; see, for example, [17]. But note that we consider arbitrary designs of $\Xi$ and not only designs with finite support. For details, see [15].

THEOREM A.1. *A design $\xi \in \Upsilon_v[r]$ is $\phi_p$-optimal $[p \in (-\infty, 1)]$ for $K^\top \boldsymbol{\theta}$ in $\Upsilon_v[r]$ if and only if the inequality*

(A.1) $$\forall y \in \mathcal{E} \ \forall z \in S : f(y)^\top N f(y) \leq f(z)^\top N f(z)$$

TABLE 1
$\alpha_0^{(k)}$ *of Example* 3.4

| $k$ | 2 | 3 | 4 | 5 | 6 | 7 | 8 |
|---|---|---|---|---|---|---|---|
| $\alpha_0^{(k)}$ | 1 | 0.75 | 0.8333 | 0.7721 | 0.7980 | 0.7755 | 0.7882 |

holds, where $S := \{z \in \mathcal{E} | \forall \epsilon > 0 : (\xi - rv \cdot \tilde{\lambda})((z - \epsilon, z + \epsilon) \cap \mathcal{E}) > 0\}$ and $N = M(\xi)^{-1} K (K^\top M(\xi)^{-1} K)^{-p-1} K^\top M(\xi)^{-1}$.

PROOF OF THEOREM 3.1. We consider Theorem A.1 with $\mathcal{E} = [-1, 1]$, $p = -1$ and $K = \mathbf{e}_k$. Then we have $N = M(\xi)^{-1} \mathbf{e}_k \mathbf{e}_k^\top M(\xi)^{-1}$ and $f(x)^\top N f(x) = (\mathbf{e}_k^\top M(\xi)^{-1} f(x))^2$ for all $x \in [-1, 1]$. Note that $\mathbf{e}_k^\top M(\xi)^{-1} f(y)$ is a polynomial of degree $k-1$ since $M(\xi)^{-1}$ is positive definite. Hence, this polynomial has at most $k$ extremal points in $[-1, 1]$. Thus, the assertion follows. $\square$

In the proof of Theorem 3.2 we use the following lemma.

LEMMA A.2. Let $k \geq 2$, let $F = ((x_{i-1})^{\ell-1})_{\ell,i=1}^k \in \mathbb{R}^{k \times k}$, where $x_i = \cos(\frac{(k-1-i)\pi}{k-1})$, $i = 0, \ldots, k-1$, and let $S = \mathrm{diag}(\ldots, -1, 1, -1, 1)$. Then we have

$$\forall x \in [-1, 1] : \mathbf{1}_k^\top S F^{-1} f(x) = T_{k-1}(x).$$

PROOF. The function $x \mapsto \mathbf{1}_k^\top S F^{-1} f(x)$ is a polynomial of degree at most $k - 1$. We have

$$\mathbf{1}_k^\top S F^{-1} f(x_i) = \mathbf{1}_k^\top S \mathbf{e}_{i+1} = (-1)^{k-1-i} = \cos((k-1-i)\pi)$$
$$= \cos((k-1)\arccos(x_i)) = T_{k-1}(x_i), \quad i = 0, \ldots, k-1.$$

Hence the assertion follows. $\square$

PROOF OF THEOREM 3.2. First note that $r \in [0, \alpha_0^{(k)}]$ is a natural condition because otherwise $\xi$ is not a measure any more. Let $x_i = \cos(\frac{(k-1-i)\pi}{k-1})$, $i = 0, \ldots, k-1$, be defined as in Lemma A.2 and let $p_i^* := \frac{\nu_i}{k-1}$. Then $\xi^* := \sum_{i=0}^{k-1} p_i^* \delta_{x_i}$ is an $\mathbf{e}_k$-optimal design in $\Xi$. Hence,

$$\forall y \in [-1, 1] \; \forall i \in \{0, 1, \ldots, k-1\} : (\mathbf{e}_k^\top M(\xi^*)^{-1} f(y))^2 \leq (\mathbf{e}_k^\top M(\xi^*)^{-1} f(x_i))^2.$$

Thus, to show Theorem 3.2 it is sufficient by Theorem 3.1 to show, for some $\gamma \in \mathbb{R}$,

(A.2) $$M(\xi^*)^{-1} \mathbf{e}_k = \gamma M(\xi)^{-1} \mathbf{e}_k.$$

Let $F, S = \mathrm{diag}(\ldots, -1, 1, -1, 1)$ be defined as in Lemma A.2, let $R = \int_{-1}^1 f f^\top \times v \, d\tilde{\lambda}$, and let $\tilde{F}$ and $\tilde{R} \in \mathbb{R}^{(k-1) \times k}$ be the first $k - 1$ rows of $F$ and $R$, respectively. Further, let $\mathbf{p}^* := (p_0^*, \ldots, p_{k-1}^*)^\top$, $\mathbf{q} := (q_0, \ldots, q_{k-1})^\top$, where $q_i$ is defined in Theorem 3.2, and let $P^* := \mathrm{diag}(p_0^*, \ldots, p_{k-1}^*)$, $Q := \mathrm{diag}(q_0, \ldots, q_{k-1})$. Note that $M(\xi^*) = F P^* F^\top$, $M(\xi) = M(\xi^*) + r(R - F Q F^\top)$. Hence (A.2) is equivalent to

$$(R - F Q F^\top) F^{\top -1} P^{*-1} F^{-1} \mathbf{e}_k = \gamma \mathbf{e}_k \qquad \text{for some } \gamma \in \mathbb{R}.$$

Thus, to prove the assertion it is sufficient to show

$$\mathbf{0}_{k-1} = (\tilde{R} F^{\top-1} P^{*-1} F^{-1} - \tilde{F} Q P^{*-1} F^{-1}) \mathbf{e}_k$$

$$(\text{A.3}) \qquad = (\tilde{R} F^{\top^{-1}} S P^{*-1} S F^{-1} - \tilde{F} S Q P^{*-1} S F^{-1}) \mathbf{e}_k$$

$$= \tilde{R} F^{\top^{-1}} S \mathbf{1}_k - \tilde{F} S \mathbf{q},$$

where the last equation holds since Elfving's theorem ($\mathbf{e}_k$-optimality of $\xi^*$ in $\Xi$) gives $\frac{\mathbf{e}_k}{\rho(\mathbf{e}_k)} = F S \mathbf{p}^*$. We define the following functions and matrices using Chebyshev polynomials:

$$h^{(1)}(x) = (T_0(x), \ldots, T_{k-2}(x))^\top,$$
$$h^{(2)}(x) = (T_{k-1}(x), \ldots, T_{2k-3}(x))^\top,$$
$$h(x) = (h^{(1)}(x)^\top, h^{(2)}(x)^\top)^\top,$$
$$H^{(j)} = (h^{(j)}(x_0), \ldots, h^{(j)}(x_{k-1})), \qquad j = 1, 2,$$
$$H = (h(x_0), \ldots, h(x_{k-1})).$$

With these definitions we can write

$$\mathbf{q} = \frac{1}{2(k-1)} \operatorname{diag}(1/2, 1, 1, \ldots, 1, 1/2) H^\top \int_{-1}^1 h(x) v(x) \, dx.$$

Since $(-1)^{k-1-i} T_j(x_i) = T_{j+k-1}(x_i)$, we get by using problems 1.5.28 and 1.1.3 in [19],

$$H^{(1)} S \mathbf{q} = \frac{1}{2(k-1)} H^{(2)} \operatorname{diag}(1/2, 1, 1, \ldots, 1, 1/2) H^\top \int_{-1}^1 h(x) v(x) \, dx$$

$$(\text{A.4}) \qquad = \begin{pmatrix} 0 & 0 & \cdots & 0 & 1/2 & 0 & \cdots & 0 \\ 0 & 0 & & 1/4 & 0 & 1/4 & & 0 \\ \vdots & & & & & & \ddots & \\ 0 & 1/4 & \cdot\cdot & 0 & & \cdots & & 0 & 1/4 \end{pmatrix} \int_{-1}^1 h(x) v(x) \, dx$$

$$= \frac{1}{2} \int_{-1}^1 h^{(1)}(x) T_{k-1}(x) v(x) \, dx.$$

Next we multiply (A.4) on the left side by the matrix which changes the basis from $T_0, \ldots, T_{k-2}$ to $x^0, \ldots, x^{k-2}$ and use Lemma A.2. Then we obtain $\tilde{F} S \mathbf{q} = \tilde{R} F^{\top^{-1}} S \mathbf{1}_k$, implying that (A.3) is true. $\square$

**Acknowledgments.** The first version of this paper was written while both authors had positions at the Institute of Mathematical Stochastics, University of Karlsruhe, Germany. We thank an Associate Editor and two referees for valuable hints and suggestions.

FACULTY OF MATHEMATICS AND GEOGRAPHY
CATHOLIC UNIVERSITY OF EICHSTÄTT–INGOLSTADT
D-85071 EICHSTÄTT
GERMANY
E-MAIL: wolfgang.bischoff@ku-eichstaett.de

CLINICAL INFORMATION SCIENCE
ASTRAZENECA
S-15185 SÖDERTÄLJE
SWEDEN